\input amstex
\loadmsbm
\loadbold
\documentstyle{amsppt}

\magnification=1200
\define\Spec{\operatorname{Spec}}

\define\OK{\Cal O_K}
\define\E{\overline{E}}
\define\T{\overline{T}}
\redefine\P{\Bbb P}
\define\F{\Bbb F}
\define\n{\underline{n}}
\define\U{\frak U}
\redefine\*{\ast}
\define\<{\langle}
\define\>{\rangle}

\topmatter
\title{Heights and Geometric Invariant Theory}
\endtitle
\date 
January 20th, 1997\enddate

\author Carlo Gasbarri\endauthor
\address Universit\'e de Rennes 1, D\'ept. de Math\'ematiques,
 Campus de Beaulieau, 35000 Rennes (F)\endaddress
\email gasbarri\@forestiere.univ--rennes1.fr\endemail

\abstract Let $K$ be a number field, $\OK$ be its ring of
integers. We introduce the notion of compactified representation
of $GL_N(\OK)$ and, we see how to associate to a hermitian 
vector bundle $\E$ over $\Spec(\OK)$ and a compactified representation
$\T$, a hermitian tensor bundle $\E_T$. We can prove then 
that there exists a lower bound
for the heights of points $x\in\P(\E_T)$ with $SL_N(K)$--semistable
generic fibre
in terms of the degree of $\E$ and some universal constants depending
only on the compactified representation. We give then three 
applications: a universal lower bound for general flag varieties,
an application to the adjoint representation of $SL_N(K)$ and
a construction of a height on the moduli space of semistable
vector bundles over algebraic curves.
\endabstract
\endtopmatter

\document

{\bf \S 1 Introduction}

\

In the papers [Bo1], [Bo2], [Bu], [So], [Zh1], [Zh2], the different
authors (J.B. Bost, J.F Burnol, C. Soul\'e, S. Zhang) have shown the 
interesting relations between height theory and geometric 
invariant theory. In particular J.B. Bost [Bo1] and S. Zhang [Zh2]
shown that, if $X\subset\Bbb P^N(\overline{\Bbb Q})$ is a closed 
variety which has $SL(N+1)$--semistable Chow point then the
height of $X$ cannot be too small.

Here we continue their program; we consider a (quite) arbitrary
 linear action of $GL_N(K)$ on a projective space $\Bbb P^M$
and we study the height of $SL_N(K)$--semistable or unstable
points under this action.

Let $K$ be a number field and $\Cal O_K$ be its ring of integers.

Let $T:GL_N(\OK)\to GL(W)$ be a linear representation and
$\overline{E}\to \Spec(\OK)$ be an hermitian vector bundle of rank
$N$. As in classical algebraic geometry (the function field case), 
 from $\overline{E}$ and $T$ we would like to construct an associated
hermitian tensor bundle $\E_T$.  The problem is:
what is the metric on $E_T$?

We solve this problem by giving a definition which fits very well
in the Arakelov geometry: the notion of 
{\it compactified representation} $\overline T$ of $GL_N(\OK)$
(cfr. chapter 3). 
So, given a compactified representation $\T$ and an hermitian vector
bundle $\E$ we can uniquely construct an associated hermitian
tensor bundle $\E_T$.

Let $a\in\Bbb Z$, we will say that a representation 
$T:GL_N(K)\to GL(W)$ is homogeneous of degree $a$ if, for all
$t\in K^{\*}$ and all $w\in W$, $T(t\cdot Id_N)(w)=t^a\cdot w$.

We can now state the main Theorem of this paper:

\proclaim{Theorem 1} Let 
$\T =(T:GL_N(\OK )\to GL(W);\, \lambda_{\sigma})_{\sigma}$
be a homogeneous compactified representation of degree $a$; then
there exists a explicitly computable constant $C=C(\T)$ such that,
for all hermitian vector bundle $\E\to\Spec(\OK)$ and
all point $p\in\Bbb P(\E_T)(\OK)$ which  has restriction
to the generic fibre $p_K$ which is $SL_N(K)$--semistable, we have
$$h(p)\geq{{a}\over{[K:\Bbb Q]}}\cdot {{\widehat{\deg }(\E)}\over {N}}
+C.$$
$h(p)$ is the height of $p$ in $\Bbb P(\E_T)$.

Conversely: let $x\in\Bbb (W_K)$ be an $SL_N(K)$--unstable point,
then there exist a sequence of vector bundles of rank $N$,
$\E_n\to Spec(\OK)$, such that, if $x_{\E_n}\in\Bbb P(\E_n)$ is
the point corresponding to $x$, then
$$\lim_{n\to\infty}\left (
h(x_{\E_n})-
{{a}\over{[K:\Bbb Q]}}\cdot {{\widehat{\deg }(\E_n)}\over {N}}\right )=
-\infty.$$\endproclaim

We can prove also an analogue statement for arbitrary compactified
representations (Theorem 2).

We propose then three applications of Theorem 1:

-- First we can give a lower bound for the height of the flag varieties
in term of the degree of the hermitian vector bundle:

\proclaim{Proposition 5.1} Let $N$ be a positive integer and
$\underline{n}=(n_1;\dots ;n_k)$ be a partition of $N$; There exist
two universal constants $A(\underline{n})$ and $B(\underline{n})$
such that, if $\E\to\Spec (\OK )$ is an hermitian vector bundle
of rank $N$ and $\Bbb F(\underline{n})(\E )$ is the flag variety
associated to $\E$ and the partition $\underline{n}$, then
$$h(\Bbb F(\underline{n})(\E ))\geq A(\underline{n})
{{\widehat{\deg (\E)}}\over{N}}+B(\underline{n}).$$
\endproclaim

(See chapter 5, paragraph 1 for more details).

-- As second application we can consider heights on the projective
space associated to the Lie algebra $\frak{sl}(N)$, 
$\Bbb P(\frak{sl}(N))$. Using a caracterisation (due to Mumford)
of semistable points on $\Bbb P(\frak{sl}(N))$ and Theorem 1,
we can give a universal lower bound for the height of those
points $x\in\Bbb P(\frak{sl}(N))$ such that $ad(x)$ is
not nilpotent. (See chapter 5 paragraph 2).

-- The third application is a bit different from the others.
Using the work of J.F. Burnol [Bu] and S. Zhang [Za1], we can
construct a (non canonical) height on the moduli space of
semistable vector bundles over an algebraic curve.
This height seems quite computable and it is very natural in terms
of the geometric invariant theoretic construction of
such a moduli space. This third application may also be seen as
an appendix to another paper on the argument [Ga]. In that paper we 
constructed a ''canonical'' height under some hypothesis on the 
curve (good reduction everywhere) and the degree and the rank
of the bundles. That height has a very interesting interpretation
in terms of the Arakelov geometry of the bundles. The height
constructed here has (a priori) no Arakelovian interpretation,
but we do not make any hypothesis on the curve, the degree 
and the rank. It might be interesting to study the relations
between the two heights. 

We would like to thank the Equipe de G\'eom\'etrie Alg\'ebrique
of the Rennes 1 University for its warm and friendly hospitality and 
J.B. Bost, B. Edixhoven, L. Moret--Bailly, C. Soul\'e, L. Szpiro,
 S. Zhang, for helpful discussions.

This research is partially supported by a
 EEC--HCM PostDoctoral Fellowship.

\

{\bf \S 2 Quick review of heights theory of projective varieties
and Geometric Invariant Theory}

\

{\it -- Heights theory}

\

Let $K$ be a number field, $\OK$ be its ring of integers
and $S_{\infty}$ the set of infinite places of $K$.

An hermitian vector bundle of rank $N$, $\E$ over
$\Spec (\OK)$ is a couple 
$\E=(E;\<\cdot ;\cdot\>_{\sigma})_{\sigma\in S_{\infty}}$ where:

-- $E$ is a projective $\OK$--module of rank $N$ (or, equivalently,
a vector bundle over the arithmetic curve $\Spec (\OK)$).

-- For each $\sigma\in S_{\infty}$, the $\Bbb C$--vector space
$E\otimes_{\sigma}\Bbb C=E_{\sigma}$ is endowed with a hermitian
metric $\<\cdot ;\cdot\>_{\sigma}$ (if
$\sigma$ is the complex conjugate of $\tau$ then the metric
$\<\cdot ;\cdot\>_{\sigma}$ is the complex conjugate
of $\<\cdot ;\cdot\>_{\tau}$).

We can associate to $\E$ an arithmetic variety $X=\Bbb P(\E)$.
More precisely, in the language of Gillet and Soul\'e arithmetic
intersection theory, $X$ is an ''Arakelov Variety''; indeed, for
each $\sigma\in S_{\infty}$, 
$X_{\sigma}=X\times_{\sigma}\Spec (\Bbb C)$ is endowed with an 
hermitian metric: the Fubini--Study metric. Moreover $X$ is
equipped with a canonical universal hermitian quotient line bundle
$\overline{L}=\overline{\Cal O_X(1)}$.

In this situation (actually in a more general one) Gillet
and Soul\'e defined an arithmetic intersection theory with values in
$\Bbb R$ (cfr. [GS]); actually we have:

-- for each integer $i\in [O;N]$, an arithmetic
Chow group $\widehat{CH^i(X)}$;

-- an intersection product: 
$(\cdot ;\cdot):\widehat{CH^i(X)}\otimes\widehat{CH^j(X)}\to
\widehat{CH^{i+j}(X)}$;

-- an arithmetic degree map: 
$\widehat{\deg }:\widehat{CH^N(X)}\to \Bbb R$;

-- for each cycle of dimension $p$, $Z\in Z_p(X)$, an arithmetic 
fundamental class $[Z]\in \widehat{CH^{N-p}(X)}$ (this class depends
on the chosen metric on $E$).

By all this machinery we can define a height function on cycles of
arbitrary dimension on $X$:
let $Z\in Z_p(X)$, we define the height of $Z$ by the formula
$$h(Z)={{1}\over{[K:\Bbb Q]}}\widehat{\deg}(\overline{L}^d;[Z]).\tag 1$$
For an extended account on height theory see [BoGS], [F].

\remark{Remark} Suppose $Z_K\in Z_{p-1}(X_K)$ is a cycle of dimension
$p-1$ on the generic fibre $X_K$ of $X$; then we define its degree
 by the formula $\deg(Z_K)=\deg(L_K^{p-1};[Z_K])$ (where $[Z_K]$
is the geometric fundamental class). comparison between this formula and
formula (1), shows that we can see the height as the arithmetic 
analogue of the geometric degree. And, if the last one measure the
geometric complexity, the former (by analogy) measure the arithmetic
complexity of a cycle in $X$.\endremark

{\it Example:} Let $p\in X(\OK)$; then $p$ correspond to  projective
module $M$ of rank one over $\OK$ with a surjection
$$E\longrightarrow M\to 0$$
if we put over $M$ the quotient metric induced by $\E$, we have the
formula
$$h(p)={{1}\over{[K:\Bbb Q]}}\widehat{\deg}(M)$$
where $\widehat{\deg}(M)$ is the Arakelov degree defined by the 
formula
$$\widehat{\deg}(M)=\log (M\big/_{m\OK})
-\sum_{\sigma}\log\Vert m\Vert_{\sigma}$$
where $m\in M\setminus\{ 0\}$ and 
$\Vert m\Vert_{\sigma}=\< m;m\>_{\sigma}^{1/2}$ (cfr. [Sz]).

\

{\it -- Geometric Invariant Theory}

\

Let $G$ be the reductive group $GL_N(K)$ or $SL_N(K)$.

Let
$$T:G\longrightarrow GL(W)$$
be an algebraic linear representation of $G$, where $W$ is a finite
dimensional $K$--vector space.

\definition{Definition} A point $v\in W$ is said
{\it semistable} with respect to $T$ if it verify one of 
these two equivalent conditions:

a) the Zariski closure of the orbit of $v$, 
$\overline {T(G)(v)}$ does not contains the zero element;

b) there exists a $G$--invariant homogeneous polynomial of
 positive degree $f$ such that $f(v)\neq 0$ 
($f\in Sym^{\*}(W^{\*})^G$).\enddefinition

If $v$ is not semistable, it is said unstable.

By linearity of the action $T$, the group $G$ acts on the projective
space $\Bbb P(W)$. A point $p\in\Bbb P(W)$ is said {\it semistable}
if there exists a vector $v\in W$ over $p$ which is semistable
(and then all the other vectors $v'$ over $p$ are semistable too).
Otherwise $p$ is said unstable.

The following facts hold:

-- The set of semistable points of $W$ (resp. $\Bbb P(W)$) is a Zariski
open set $W^{ss}$ (resp. $\Bbb P(W)^{ss}$). This open set is 
invariant under
base change: if $K'$ is an extension of $K$ and 
$T':G\otimes_KK'\to GL(W\otimes_KK')$ is the extended representation, 
then
$$W^{ss}\otimes_KK'=(W\otimes_KK')^{ss}\text{   and   }
\Bbb P(W)^{ss}\otimes_KK'=(\Bbb P(W\otimes_KK')^{ss}.$$

-- The $K$--algebra $Sym^{\*}(W^{\*})^G$ of the $G$--invariant 
polynomials over $W$, is a finitely generated $K$--algebra
 (and, as before, commute to base change).

\remark{Remark} The theory can be extended to arbitrary reductive
algebraic group; but we will not use that generality here.
\endremark

For more details on Geometric Invariant Theory see [GIT] or [Se].

\newpage

{\bf \S 3 Representations of $GL(N)$ and hermitian vector bundles}

\

Let $K$ be a field and $X$ be a variety over $K$; let $E\to X$
be a vector bundle of rank $N$ and $T:GL_N(K)\to GL(W)$ be an
algebraic linear representation defined over $K$. It is well known
that we can construct a new vector bundle $E_T$ over $X$ by using $E$
and $T$. It is  the so--called ''tensor bundle 
associated to $E$ and $T$''.

Roughly speaking $E_T$ is constructed in the following way:
$E$ can be described by an open covering $\Cal U=\{U_i\}_{i\in I}$
and algebraic functions $g_{ij}:U_i\cap U_j\to GL_N(K)$ 
(called transition functions) which satisfy the well known cocycle
conditions; the $E_T$ is described by the same open covering $\Cal U$
and as transition functions the functions $h_{ij}=T(g_{ij})$.

{\it Example} If $V$ is the standard representation of $GL_N(K)$
and $T=Sym^n(V)$ (resp. $T=\bigwedge^i(V)$) then
$E_T=Sym^n(E)$ (resp. $E_T=\bigwedge^i(E)$).

We would like to have an analogue construction over arithmetic curves.

Given a linear algebraic representation 
$T:GL_N(\OK)\to GL(W)$ and a vector bundle
$E\to \Spec(\OK)$ it is not difficult to construct a new vector
bundle $E_T\to\Spec(\OK )$ which is the  analogue of the tensor bundle 
described before; actually the same construction works.

But if $\E\to\Spec (\OK)$ is an {\it hermitian} vector bundle, can
$E_T$  be endowed with a  natural structure of {\it hermitian} 
vector bundle ? The answer is ''not in general'' !

Given an hermitian vector space $V$ (over $\Bbb C$) of dimension $N$
and a representation $T:GL_N(\Bbb C)\to GL(W)$, there are plenty
of ''natural'' metrics which depend in a ''natural way'' on the metric
on $V$.

{\it Example} Let $V$ be an hermitian space and $W=Sym^n(V)$.
Over $W$ there are at least two natural hermitian metrics:

-- There is a canonical surjective map $V^{\otimes n}\to Sym^n(V)$, 
which induces on $W$ the quotient metric of the tensor product
metric.

-- We can see W as $H^0(\Bbb P(V);\Cal O_{\Bbb P}(n))$, and over this
last one there is the $L^2$ metric induced by the Fubini--Study, 
metric on $\Bbb P(V)$ (cfr. also [BoGS] page 984).

Never the less we have the following Proposition which can be useful 
to solve the problem when the representation is irreducible:

\proclaim{Proposition 3.1} Let $T:GL_N(\Bbb C)\to GL(W)$ be
an irreducible representation, let $U(N)\subset GL_N(\Bbb C)$
be the unitary group; then, up to a constant, there exist a unique
$U(N)$--invariant hermitian metric over $W$.\endproclaim

\demo{Proof} The restriction $T\vert_{U(N)}:U(N)\to GL(W)$is still
irreducible. Let $h$ and $g$ be two $U(N)$--invariant metrics on $W$.
By Riesz representation Theorem, there exist a linear map
$A:W\to W$ such that, for all $x$ and $y$ in $W$,
$h(x;A(y))=g(x;y)$. It is not difficult to see that $A$ is a 
$U(N)$--equivariant map. By Schur Lemma there is a constant 
$\lambda$ (which must be real positive) such that
$A=\lambda\cdot Id_W$.\enddemo

This Proposition can be easily generalized to representations on which
every irreducible representation appears with multiplicity
zero or one (but the constant will be a $t$--tuple). But the following 
example shows that we cannot hope better than this.

{\it Example} Let $\Bbb C$ be the trivial representation of
$GL_N(\Bbb C)$; then every hermitian metric on 
$\Bbb C\oplus \Bbb C$ is $U(N)$--invariant.

So we propose the following definition.

We recall that if $V$ is a finite dimensional $\Bbb C$--vector 
space, the set $\Cal H(V)$ of the hermitian metrics on $V$ can be seen
as a homogeneous space in the following way: $GL(V)$ acts on 
$\Cal H(V)$: if $h\in\Cal H(V)$ and $g\in GL(V)$, then
$\alpha (g)(h)= ^tgh\overline{g}$.
The action $\alpha$ is transitive and if $h\in\Cal H(V)$, the isotropy
group of $h$ is the unitary group $U(h)$. So
$$\Cal H(V)=GL(V)\bigg/_{U(h)}.$$

Let $T:GL_N(\Bbb C)\to GL(W)$ be a linear representation and let
$h\in \Cal H(W)$ be a $U(N)$--invariant hermitian metric. Then,
by using $T$ and $h$, we can define a map from
$\Cal H(\Bbb C^N)$ to $\Cal H(W)$: since $h$ is $U(N)$--invariant,
$T(U(N))\subset U(h)$, so $T$ defines a map from
 $GL_N(\Bbb C)\big/_{U(N)}$ to $GL(W)\big/_{U(h)}$, so a map
$$\lambda_{h;T}:\Cal H(\Bbb C^N)\longrightarrow \Cal H(W).$$

Let $V$ be a hermitian vector space of dimension $N$; by fixing
an hortonormal basis, it defines a point 
$a_V\in\Cal H(\Bbb C^N)$ and a compact subgroup 
$U(V)\subset GL_N(\Bbb C)$ (the point $a_V$ and the unitary group 
$U(V)$ do not depend on the
chosen basis); 
so, by the map $\lambda_{h;T}$ it defines
an hermitian metric $\lambda_{h;T}(a_V)$ on $W$. The hermitian
metric $\lambda_{h;T}(a_V)$ is $U(V)$--invariant.

\definition{Definition}: Let $K$ be a number field, $\OK$ 
its ring of integers and $S_{\infty}$ the set of embedding of $K$ in 
$\Bbb C$. A {\it compactified representation} $\overline T$
of $GL_N(\OK )$ is a couple 
$\overline{T}=(T;h_{\sigma})_{\sigma\in S_{\infty}}$, where

-- $T:GL_N(\OK)\to GL (W)$ is a linear representation and $W$ is
a free $\OK$--module.

-- For each $\sigma\in S_{\infty}$, $h_{\sigma}$ is a 
$U(N)$--invariant metric.\enddefinition

Let $\E$ be an hermitian vector bundle over $\Spec(\OK)$ and
$\overline{T}=(T:GL_N(\OK)\to GL(W); 
h_{\sigma})_{\sigma\in S_{\infty}}$
be a compactified representation. From these data we can
construct a new hermitian vector bundle  $\E_T$ over $\Spec(\OK)$ as
follows:

-- If $\Cal U=\{ U_i\}_{i\in I}$ is  covering of
$\Spec(\OK)$ by affine open sets ($U_i=\Spec (A_i)$; 
$U_i\cap U_j=\Spec (A_{ij})$) where $E$ trivializes, and
$g_{ij}:U_i\cap U_j\to GL_N(A_{ij})$ are the transition functions
defining $E$; then $E_T$ is defined by the same open covering
and as transition functions 
$h_{ij}=T(g_{ij}):U_i\cap U_j\to GL(W\vert_{U_i\cap U_j})$.

-- For each $\sigma\in S_{\infty}$ the hermitian metric 
on $(\E_T)_{\sigma}\simeq W_{\sigma}$ is 
$\lambda_{h_{\sigma};T}(E_{\sigma})$.

\definition{Definition} The hermitian vector bundle $\E_T$ will
be called {\it the hermitian tensor bundle associated to $\E$
and $\overline{T}$}.\enddefinition

{\it Examples} -- If $V$ is the standard representation of $GL_N(\OK)$
and $h_{\sigma}$  "is the identity'', then $\E_T\simeq \E$.
 
-- If $T:GL_N(\OK)\to GL(W)$ is an irreducible representation, 
then we can prove that there exists a positive integer $d$, a integer
$k$ and a $GL_N(\OK)$--equivariant surjective map
$$\varphi:(V^{\otimes d})\otimes (\det(V)^{\otimes k})
\longrightarrow W\to 0$$
(cfr. the Weyl construction, [FH] chapter 6). So, if we fix $\varphi$
we can define $h_{\sigma}$ as the quotient metric induced by
$\varphi$; then for every hermitian vector bundle of rank $N$, $\E$,
 the metric on $\E_T$ is the quotient metric given
 by the the induced map
$$\varphi_{\E}:(\E^{\otimes d})\otimes (\det(\E)^{\otimes k})
\longrightarrow \E_T\to 0.$$
By Proposition 3.1 every other ''compactification'' of $T$ is
obtained by multiplying by a positive real number 
such a compactification.

\proclaim{Corollary 3.2} Let $T$ be a irreducible representation
of $GL_N(\OK)$ and $\E$ be a hermitian vector bundle of rank $N$; 
then there exists a unique metric on $E_T$ which is 
$U(E_{\sigma})$--invariant for every $\sigma\in S_{\infty}$ and the 
isomorphism
$$\det(E_T)\simeq (\det(\E )^{{{a\cdot rk(W)}\over{N}}}$$
is an isometry ($a$ is the integer such that, if $t\in \Bbb G_m(K)$
and $v\in W_K$, then $T(t\cdot Id_N)(v)=t^{a}\cdot v$;
 it depends only on $T$, cfr. below).\endproclaim

\remark{Remark} Let $T:GL_N(\Bbb C)\to GL(W)$ be a linear
 representation, and $h_1$ and $h_2$ be two $U(N)$--invariant 
metrics; then there exists two positive constants
$C_1=C_1(h_1;h_2)$ and $C_2=C_2(h_1;h_2)$ such that {\it for every}
$N$--dimensional hermitian vector space $V$ and element $x\in W$
$$C_1\Vert x\Vert_{\lambda_{h_1;T}(a_V)}\leq
\Vert x\Vert_{\lambda_{h_2;T}(a_V)}
\leq C_2\Vert x\Vert_{\lambda_{h_1;T}(a_V)}$$
where $\Vert x\Vert_{\lambda_{h_i;T}(a_V)}=
(\lambda_{h_i;T}(a_V)(x;x))^{1/2}$.\endremark

\

{\bf \S 4 Heights of semistable points.}

\

In this chapter, we will prove, roughly speaking, that the height of a
 semistable
point in the projective space $\Bbb P(\E_T)$ cannot be too small.

\definition{Definition} Let $K$ be a field (of characteristic zero)
and $T:GL_N(K)\to GL(W)$ be a linear representation; $T$ is said to
be {\it homogeneous of degree $a$}, where $a$ is a integer, if, for
all $t\in K^{\*}$ and $v\in W$, we have
$$T(t\cdot Id_N)(v)=t^a\cdot v.$$\enddefinition

The class of homogeneous representations is quite big, indeed:

-- every irreducible representation is homogeneous;

-- if $T$ is homogeneous (of degree $a$) then $Sym^n(T)$, 
$\bigwedge^i(T)$, $T^{\*}$ are homogeneous (of degree $na$, $ia$, $-a$
respectively);

-- every subrepresentation of a homogeneous representation is 
homogeneous (so, using duality, every quotient representation of
a homogeneous representation is homogeneous);

-- if $T$ and $T'$ are homogeneous representations, then 
$T\otimes T'$ is homogeneous;

The direct sum of homogeneous representation of different degree
{\it is not} homogeneous (and this is, more or less, the only method
to construct non homogeneous representations).

Let $K$ be our number field and $\OK$ be its ring of integers.

We can now state the principal Theorem of this paper:

\proclaim{Theorem 1} Let 
$\T =(T:GL_N(\OK)\to GL(W); \, h_{\sigma})_{\sigma\in S_{\infty}}$ be
a compactified homogeneous representation of degree $a$; then there 
exist  constant $C=C(\T)$ such that, for every hermitian vector
bundle $\E\to\Spec(\OK)$ of rank $N$ and for every point
$p\in\Bbb P(\E_T)(\OK)$ such that, the restriction $p_K$ of $p$ to the
generic fibre is $SL_N(K)$--semistable, we have
$$h(p)\geq{{a}\over{[K:\Bbb Q]}}\cdot{{\widehat{\deg}(\E)}\over{N}}
+C.$$

Conversely: let $x\in\Bbb P(W_K)$.
Let 
$$\Cal V=\left\{ \E'\text{ herm. vect. bun. of rk $N$ over
 } \Spec(\Cal O_{K})\right\},$$ 
then, for every $\E'\in\Cal V$, $x$ defines a point
$x_{\E'}$ in $\Bbb P(\E'_T)$. If $x$ is a $SL_N(K)$--unstable point
then
$$\inf_{\Cal V}\left\{
h(x_{\E'})-{{a}\over{[K:\Bbb Q]}}\cdot{{\widehat{\deg}(\E')}\over{N}}
\right\}=-\infty.$$\endproclaim

\remark{Remark} The proof of this Theorem is very similar to the
proof of Proposition 2.1 in [Bo1] and Proposition 4.2 in [Zh2].
\endremark
\demo{Proof} Let $S=Sym^{\*}(W_K)$ be the symmetric algebra of
$W$; $SL_N(K)$ acts linearly over $S$;
let $S^{SL_N}$ be the subalgebra of the $SL_N(K)$--invariant elements
of $S$. By the first fundamental Theorem in Geometric Invariant Theory,
$S^{SL_N}$ is a $K$--algebra of finite type. Let 
$\left\{ P_1;\dots;P_M\right\}$ be a basis of $S^{SL_N}$ over $K$.
By clearing the denominators, we can suppose that 
$P_i\in S(W)$ (the $\OK$--symmetric algebra of $W$); 
moreover we can suppose that the $P_i$'s are homogeneous.

Let $\sigma\in S_{\infty}$ and $h_{\sigma}$ the  hermitian metric
on $W_{\sigma}^{\*}$  (the dual of $W_{\sigma}$) given by the 
dual compactification; we define then
$$\vert\Vert P_i\Vert\vert_{\sigma}=\sup_{v\in W_{\sigma}^{\*}}
{{\Vert P_i(v)\Vert}\over{\Vert v\Vert_{h_{\sigma}}^{\deg (P_i)}}}$$
and
$$B_{\sigma}=
\max_{1\leq i\leq M}\left\{ 
\vert\Vert P_i\Vert\vert_{\sigma}^{1/\deg(P_i)}\right\}.$$
We remark that, for all $v\in W_{\sigma}^{\*}$, $P_i(v)\in\Bbb C$, 
then we can speak about its norm.

Let $\E\to\Spec(\OK)$ be an hermitian vector bundle of rank
$N$.

Let $p\in \Bbb P(\E_T)$ a point with restriction to the generic
fibre $p_K$ which is $SL_N(K)$--semistable. By functoriality
$p$ defines  a line bundle $M$ over $\Spec(\OK)$ with a surjection
$$\E_T\longrightarrow M\to 0;\tag 3$$
as before we put on $M$ the quotient metric.

By dualizing (3) and tensorizing by $M$, we get an isometric embedding
$$0\to\overline{\Cal O}\buildrel{\varphi}\over\longrightarrow 
(\E_T)^{\*}\otimes M\tag 4$$
where $\overline{\Cal O}$ is the trivial line bundle with trivial
metrics ($\OK;\Vert 1\Vert_{\sigma}=1$).

Since $p_K$ is $SL_N(K)$--semistable, there exists a 
$SL_N(K)$--invariant polynomial of positive degree $P$ such that
$P(\varphi_K)\neq 0$. We can suppose that $P$ is one of the $P_i$'s.

We remark that, for every $\OK$--algebra $A$, by tensor product,
the representation $T$ induces a representation
$$T_A:GL_N(A)\longrightarrow GL(W\otimes_{\OK}A).$$
Since $P$ is homogeneous (of degree, say $D$) and $SL_N$--invariant;
and since the only characters of $GL_N$ are the tensor powers of the
determinant, for every $\OK$--algebra A, every $v\in W\otimes_{\OK}A$
and $g\in GL_N(A)$ we have
$$P(T_A(g)(v))=(\det(g))^{{aD}\over{N}}P(v).$$ 
Suppose that $\E$ and $M$ are defined by the open  affine covering
$\Cal U=\{ U_i\}_{i\in I}$ (where we can suppose that
$U_i=\Spec (A_i)$ and $U_i\cap U_j=\Spec (A_{ij})$) and transition 
functions $g_{ij}$ and $f_{ij}$ respectively. Then $\varphi$ is given
by sections
$$\varphi_i\in\Gamma(U_i;(E_T)^{\*}\otimes M)\simeq W\otimes_{\OK}A_i$$
with the relations
$$\varphi_i=T^{\*}(g_{ij})\cdot f_{ij}(\varphi_j)$$
over $U_i\cap U_j$ ($T^{\*}$ is the dual representation of $T$).

Then we see that $P(\varphi)$ defines a non zero element in
$(\bigwedge^N(E^{\*}))^{{aD}\over{N}}\otimes M^D$.

Since the norm $\Vert \varphi\Vert_{\sigma}=1$ (because of
the isometry (4)) we have that
$$\Vert P(\varphi)\Vert_{\sigma}\leq B_{\sigma}^D.$$
So
$$\eqalign{Dh(p)-{{1}\over{[K:\Bbb Q]}}\cdot Da\cdot 
{{\widehat{\deg}(\E)}\over {N}}&=
{{\widehat{\deg}(\bigwedge^N(E^{\*}))^{{aD}\over{N}}\otimes M^D}\over
{[K:\Bbb Q]}}\cr
&\geq -\sum_{\sigma}D\log B_{\sigma}\cr}.$$
So we define $C(\T)$ to be $-\sum_{\sigma}\log B_{\sigma}$.

Now we prove the converse.

Let $\Cal A=\left\{ \text {hermitian metrics } 
k=(k_{\sigma})_{\sigma\in S_{\infty}} \text{ on } 
\Cal O_{K}^{\oplus N} \text{ such that } \det (\Cal O_{K}^{\oplus N};
k)\simeq \overline{\Cal O_{K}}\right\}$. 

For each element
$k\in\Cal A$ we will denote $\E_k$ the corresponding hermitian
vector bundle over $\Spec(\Cal O_{K})$.

the set $\Cal A$ can be seen as a ''homogeneous space''
$$\prod_{\sigma\in S_{\infty}}SL_N(\Bbb C)\Bigg/_{\prod_{\sigma\in 
S_{\infty}}SU(N)}$$
where $\prod_{\sigma\in S_{\infty}}SU(N)$ acts on
$\prod_{\sigma\in S_{\infty}}SL_N(\Bbb C)$ on the left. This can be 
seen in the following way: let $e_1;\dots;e_N$ be the standard basis
of $(K)^N$: then we send the element
$(g_{\sigma})_{\sigma\in S_{\infty}}
\in\prod_{\sigma\in S_{\infty}}SL_N(\Bbb C)$ to the metric
$k=(k_{\sigma})_{\sigma\in S_{\infty}}$ having 
$\{ g_{\sigma}(e_i)\}$ as hortonormal basis.

Let $k\in\Cal A$ and $\E_k$ the corresponding vector bundle;
the associated tensor bundle $(E_k)_T$ (without metrics) is just $W$.

The action of $SL_N(\OK)$ on $\Bbb P(W)=\Bbb P((\E_k)_T)$ can
be described as a morphism
$$SL_N(\OK)\times\Bbb P(W)\buildrel{T}\over\longrightarrow\Bbb P(W)$$
which satisfies some cocycle conditions (cfr. [GIT]); let 
$L=\Cal O_{\Bbb P}(1)$ be the universal quotient bundle on
$\Bbb P(W)$ and let $p_2:SL_N(\OK)\times\Bbb P(W)\to\Bbb P(W)$ 
be the second projection. Since the action of $SL_N(\OK)$ is
linear, we have an isomorphism
$$\phi:p_2^{\*}(L)\simeq T^{\*}(L).\tag 5$$

The metric $k$ induces a metric $\Vert\cdot\Vert_k$ on L; let
$\overline L_k$ be the corresponding hermitian line bundle on
$\Bbb P(W)$. The isomorphism (5) is not an isometry in general, then
$$p_2^{\*}(\overline L_k)\simeq T^{\*}(\overline L_k)
\otimes \Cal O(\mu)$$
where $\Cal O(\mu)$ is the trivial bundle with norm
$\Vert 1\Vert_{\sigma}(g;x)=\exp (-\mu_{\sigma}(g;x))$ where
$\mu_{\sigma}:SL_N(\Bbb C)\times \Bbb P(W_{\sigma})\to \Bbb R$ is
a function.

Now, let $x\in \Bbb P(W_K)$; let $h(x)$ be its height when we see
it as a point of $\Bbb P((\overline{\Cal O^{\oplus N}}))$; let 
$g=(g_{\sigma})_{\sigma\in S{\infty}}$ be an element in $\Cal A$ and
$\E_g$ the corresponding hermitian vector bundle; let $x_{\E_g}$ be the
point $x$ when we see it as a point of
$\Bbb P((\E_g)_T)$; we have then
$$h(x_{\E_g})=h(x)+{{1}\over{[K:\Bbb Q]}}\sum_{\sigma\in S_{\infty}}
\mu_{\sigma}(g_{\sigma}^{-1};x_{\sigma}).$$
Then we conclude by using Theorem 2.2 in [Zh2].\enddemo
\remark{Remarks} 1) As we can see from the proof, the constant
$C(\T)$ is effective if we know a basis of the $SL_N(K)$--invariant
polynomials of $T^{\*}$. There are effective methods to construct
such a basis.

2) Theorem 2.2 in [Zh2], roughly speaking, says that, if we fix $x$,
the functions $\mu_{\sigma}(g;x)$ are bounded below if and only if
$x$ is semistable.\endremark

We will now shortly analyze the case of an arbitrary 
representation of $GL_N(\OK)$.

Let firstly prove an easy generalization of Proposition 3.1
\proclaim {Lemma 4.2} Let $T_i:GL_N(\Bbb C)\to GL(W_i)$, $i=1,2$
be two linear representations with no isotypic common factors;
let
$$W=W_1\oplus W_2\tag 6$$
be the direct sum representation and $h$ an $U(N)$--invariant
metric on $W$. Then the decomposition (6) is an hortogonal
decomposition with respect to $h$.\endproclaim
\demo{Proof} Let $g$ be an $U(N)$--invariant metric on $W$ for
which the (6) is an hortogonal decomposition.
Let $A:W\to W$ be a linear map such that, for all $x,y\in W$,
$h(x;y)=g(x;A(y))$. As in proposition 3.1 $A$ is an 
$U(N)$--equivariant map. The maps
$$\varphi_{ij}:W_i\buildrel{A}\over\longrightarrow W
\buildrel{pr_j}\over\longrightarrow W_j$$
($pr_j$ are the projections) are $U(N)$--equivariant maps.
But, since the $W_i$'s have no isotypic common factors,
$\varphi_{ij}=0$ if $i\neq j$. So $A(W_i)\subset W_i$ and the Lemma
is proved.\enddemo

Let $\T=(T:GL_N(\OK)\to GL(W); h_{\sigma})_{\sigma\in S_{\infty}}$
be a compactified representation; then, by the Lemma, we can
write $\T$ in a unique way as a direct sum of compactified
representation $\T=\oplus_i\T_i$, where 
$T_i:GL_(\OK)\to GL(W_i)$ are homogeneous of degree $a_i$ ($a_i\neq a_j$
if $i\neq j$).
 Let $\E\to\Spec(\OK)$ be an hermitian vector bundle of rank $N$;
again by the Lemma we can decompose the hermitian vector bundle
$E_T$ as direct sum $\E_T=\oplus_i\E_{T_i}$.

We can now state the analogue of Theorem 1 for
arbitrary representations.

\proclaim{Theorem 2} Let 
$\T=(T:GL_N(\OK)\to GL(W); h_{\sigma})_{\sigma\in S_{\infty}}$
be a compactified representation; let $\T=\oplus_i\T_i$ be the 
decomposition of $\T$ as a direct sum of homogeneous representations,
 $\deg T_i=a_i$.Let $\E\to\Spec(\OK)$ be an hermitian vector bundle
of rank $N$ and let $\E_T=\oplus \E_{T_i}$ be the hermitian
tensor bundle associated to $\E$ and $\T$. Let
$A_T(\E)=\min_i\left\{ {{a_i}\over{[K:\Bbb Q]}}\cdot
{{\widehat{\deg}(\E)}\over{N}}\right\}$. There exists a constant
$C=C(\T)$ {\rm depending only on $\T$}, such that the following holds:
let $x\in\Bbb P(\E_T)(\OK)$ be a point which verifies one
of these two properties:

a) $x\in\Bbb P(\E_{T_i})\subset\Bbb P(\E_T)$ and $x_K$ is semistable;

b) $x\not\in \Bbb P(\E_{T_i})$ for all $i$, but there exists $i$ 
such that, if $p_i:\Bbb P(\E_T)\to\Bbb (\E_{T_i})$ is the linear 
projection induced by the exact sequence
$$0\to\bigoplus_{j\neq i}\E_{T_j}\longrightarrow \E_T\longrightarrow
\E_{T_i}\to 0;$$
$p_i(x)_K$ is semistable.
Then
$$h(x)\geq A_T(\E)+C.$$
\endproclaim
\remark{Remark} If $\widehat{\deg}(\E)\geq 0$ then
$A_T(\E)=\min_i\{ a_i\}\cdot{{1}\over{[K:\Bbb Q]}}\cdot
{{\widehat{\deg}(\E)}\over{N}}$.\endremark
\demo{Proof} If $x$ verifies a) then we just apply Theorem 1 
to $x$ in $\Bbb P(\E_{T_i})$; if $x$ verifies b) we apply
Theorem 1 to $p_i(x)$ in $\Bbb P(\E_{T_i})$ and the comparison of the 
heights of projections as in [BoGS] 3.3.2.\enddemo
\remark{Remarks} a) If $x\in \Bbb P(E_{T_i})$ is unstable, we can
state an analogue of the converse of Theorem 1.

b) If $T_i$, $i=1;2$ are two representations and 
$x=x_1+x_2\in T=T_1\oplus T_2$; the if one of the $x_i$'s is 
semistable then $x$ is semistable, but in general the converse is not
true.\endremark

\

{\bf \S 5 Applications}

\

We will see three applications of Theorem 1; a lower bound for the
height of flag varieties; heights of semistable points under the
adjoint representation and a construction of a height on the
moduli space of semistable vector bundles of fixed rank and degree
over algebraic curves.

\

{\it a) Heights of flag varieties}

\

Let $N$ be a positive integer and let $\underline{n}=(n_1;\dots;n_m)$
be a partition of $N$ ($n_i\in\Bbb N_{>0}$ and $\sum n_i=N$). Let
$\E\to\Spec(\OK)$ be an hermitian vector bundle of rank $N$, and
let $\F(\n)(\E)$ be the flag variety associated to $\E$
and $\n$ (cfr. [Gr]). The variety $\F(\n)(\E)$ is the variety 
classifying flags of vector spaces of type $\n$. Namely, let
$k$ be a algebraic closed field (with a morphism $\OK\to k$)
then a closed point $x\in\F(\n)(\E)(k)$ is a flag of $k$--vector
spaces
$$\{0\}\subset F_1\subset\dots\subset F_{m-1}\subset F_m=E\otimes k$$
such that $\dim_k(F_i\big/_{F_{i-1}})=n_i$.

{\it Example} -- If $\n =(N-1;1)$ then $\F(\n)(\E)=\P(\E)$

-- If $\n =(N-p;p)$ ($p>1$) then $\F(\n)(\E)=Gr(p;\E)$,
the grassmannian of the vector bundles of rank $p$ quotient of
$\E$.

It is well known (cfr. [Gr]) that there exist a representation
$T:GL(E)\to GL(W)$ and a canonical embedding
$$i_{\n}:\F(\n)(\E)\longrightarrow \P(\E_T).$$
$T$ is of the form $\bigotimes_{i=1}^m\bigwedge^{\alpha_i}E$
(for the exact values of the $\alpha_i$ see [Gr]) then it may 
be ''quite naturally'' compactified: we put on $E_T$ the 
tensor product of the exterior product metrics. So we can speak
about the height of $\F(\n)(\E)$.

\proclaim{Proposition 5.1} There exist two universal constants
$A=A(\n)$ and $B=B(\n)$ such that
$$h(\F(\n)(\E))\geq A(\n)\widehat{\deg}(\E)+B(\n)$$
\endproclaim
\demo{Proof} The group scheme $SL(E)$ naturally 
acts on $\F(\n)(\E)$ and the embedding $i_{\n}$ is $SL(E)$--equivariant.

Let $d$ be the Krull dimension of $\F(\n)(\E)$ and $\delta$ the
degree of $\F(\n)(\E)_K$ in $\P(\E_T)_K$.

The theory of Chow forms (cfr. [Bo1] and [BoGS]) allows us to 
construct a point $\Phi_{\n}\in\P(Sym^d(E_T)^{\otimes\delta})$, the
Chow point of $\F(\n)(\E)$. $\Phi_{\n}$ is uniquely determined by 
$\F(\n)(\E)$, and conversely there exists a closed scheme
$\underline{\bold{Ch}}\in\P(Sym^d(E_T)^{\otimes\delta})$ such
that every effective cycle of dimension $d$ and degree $\delta$ in
$\P(\E_T)$ is determined by a point in
$\underline{\bold{Ch}}$.

If we put on $Sym^d(E_T)^{\otimes\delta}$ the tensor product 
of the quotient metrics induced by the surjections
$(E_T)^{d}\to Sym^d(\E_T)$ we can prove that there exists a universal
constant $C(\n)$ such that
$$h(\F(\n)(\E))\geq h(\Phi_{\n})+C(\n)\tag 7$$
cfr. [Bo] Prop. 1.3 and [BoGS] Theorem 4.3.8.

The group scheme $SL(E)$ acts on $\P(\E_T)$ fixing $\F(\n)(\E)$, so,
by functoriality of the Chow point, the point
$\Phi_{\n}$ is a fixed point of $\P(Sym^d(E_T)^{\otimes\delta})$
under the $SL(E)$ action. In particular $(\Phi_{\n})_K$ is
$SL_N$--semistable.

So applying Theorem 1 and (7) we can find two constants
$A(\n)$ and $B(\n)$ such that
$$h(\F(\n)(\E))\geq A(\n)\widehat{\deg}(\E)+B(\n).$$
\enddemo

\remark{Remarks} a) The constant $A(\n)$ is ''completely
geometric'', it depends on $\n$, $d$ and $\delta$. It do not depends
on the (natural) compactification of the 
representation used in the proof; indeed 
$A(\n)={{\prod_{i=1}^k(N-\sum_{j\leq i} n_j)d\delta}\over
{[K:\Bbb Q]\cdot N}}$. We can also explicitely compute the 
dimension $d$ et the degree $\delta$ (cfr. [Fu] example 14.6.15): 
for instance
$d= N^2-\sum_{i=1}^kn_i(\sum_{j\leq i}n_j)$; there exist also an 
explicit formula for $\delta$ but it is much more complicated
(if $\n =(N-p;p)$ then $\delta ={{1!2!\cdots (p-1)!(d-1)!}\over
{(N-p)!(N-p+1)!\cdots(N-1)!}}$).

b) The heights of grassmannians $Gr(p;\Bbb Z^N)$ have been explicitly
computed by V. Maillot [M]. Our result is less precise but it is
true for for every number field and every hermitian vector bundle
(in particular with arbitrary metric).

c) It is possible to prove Proposition 5.1 directly by using
Theorem II in [Bo1] and the results in [Ke].

d) Using the method in the proof (or more precisely, the method
used in [Bo1] Theorem II) we can find similar lower bounds for heights
of cycles on projective tensor bundles (associated to compactified 
homogeneous representations) having semistable chow points.\endremark

\

{\it b) Semistable points under the adjoint representation}

\

Let $\frak {sl}(N)$ be the Lie algebra of $SL(N)$. Let
$$Ad:SL(N)\longrightarrow GL(\frak{sl}(N))$$
be the adjoint representation; let $\overline{Ad}$ be a 
compactification of $Ad$. Let $\E$ be a hermitian vector bundle
of rank $N$; then we will denote
$\overline{\frak{sl}(\E)}$ the tensor bundle $\E_{Ad}$.

then as a direct consequence of Theorem 1 and the characterization of
semistable points of $\P(\frak{sl}(N))$ under the adjoint 
representation (cfr. [Mu] Proposition 1.15) we find

\proclaim{Proposition 5.2} There exists a constant $C\in \Bbb R$ 
depending only on the chosen compactification such that, if
$\E$ is an hermitian vector bundle of rank $N$ and
$x\in\P(\overline{\frak{sl}(\E)})$ is a point such that 
$ad(x_K)$ is not nilpotent, then
$$h(x)\geq C.$$
conversely, if $x_K\in\P(\frak{sl}(N))$ is such that $ad(x_K)$ is
nilpotent, we can find a sequence of hermitian vector bundles
$\E_n$ of rank $N$; such that, if $x_{\E_n}$ is the 
 point in $\P(\overline{\frak{sl}(\E_n)})$ defined by $x_K$, we have
$$\lim_{n\to \infty}h(x_{\E_n})=-\infty.$$
\endproclaim
\demo{Proof} It suffices to remark that the adjoint representation
is homogeneous of degree zero and $x\in\frak{sl}(N)$ is unstable
if and only if $ad(x)$ is nilpotent.\enddemo

\

{\it c) Heights on moduli space of semistable vector bundles over
algebraic curves}

\

This third application is a little bit different from the others.
Indeed, using Geometric Invariant Theory, we construct a height on the
moduli space of semistable vector bundles of fixed rank and degree
over an algebraic curve and we use Theorem 1
to give a lower bound for this height.

The height we construct is not very canonical indeed, but it has some
advantages with respect to other constructions (cfr. [Ga] for another
construction which is more canonical but needs some hypothesis on the
curve, the rank and the degree); let's quote some of them:

-- we will not make any hypothesis on the curve and on the rank and the
degree;

-- the height is strictly related with the construction of the moduli
space;

-- the height seems quite computable.

Let $X$ be a projective smooth curve of genus $g\geq 1$ over
a number field $K$.
\remark{Remark} Sometime we will silently make finite 
extensions of the base field $K$; anyway the results we are giving
are invariant under base change.\endremark

Let $\Cal E$ be a vector bundle of rank $r$ and degree $d$ over $X$;
let $\mu (\Cal E)={{d}\over{r}}$ the slope of $\Cal E$.

The vector bundle $\Cal E$ is said to be {\it semistable},
if for every subbundle $\Cal F\subset\Cal E$ we have
$$\mu (\Cal F)\leq\mu (\Cal E).$$

If we fix $r\geq 1$ and $d\in\Bbb Z$, there exists a coarse
moduli space $\U_X(r;d)$ of semistable vector bundles of
rank $r$ and degree $d$ over $X$. It is a projective variety
of dimension $r^2(g-1)+1$.

\remark{Remark} If $L$ is a line bundle of degree $n$ over $X$,
the map
$$\eqalign{\U_X(r;d)&\longrightarrow\U_X(r;d+rn)\cr
[\Cal E]&\longrightarrow [\Cal E\otimes L]\cr}$$
is an isomorphism; so if we want to study
$\U_X(r;d)$ we can suppose $d$ very big.\endremark

It is well known (cfr. [Ne] chapter 5) that, if $r$ is fixed and
$d$ is sufficiently big (how big, can be explicitly computed)
we can find a quasi projective smooth variety $R$, a 
$K$--vector space $W$, a vector bundle $\Cal U$ of rank $r$
over $X\times R$ such that the following properties are verified

-- There exist a surjective map o vector bundles over 
$X\times R$
$$W\otimes\Cal O_{X\times R}\longrightarrow \Cal U\to 0$$

-- for all $q\in R(\overline K)$, the restriction
$\Cal U_q=\Cal U\vert_{X\times\{ q\}}$ is a vector bundle of rank
$r$ and degree $d$ over $X$ and the induced map
$$W\longrightarrow H^0(X;\Cal U_q)$$
is an isomorphism;

-- If $\Cal E$ is a semistable vector bundle of rank $r$ and degree $d$
over $X$, there exists a $q\in R$ such that
$\Cal E\simeq\Cal U_q$;

-- The group $SL(W)$ acts on $R$ and $\Cal U_{q_1}\simeq\Cal U_{q_2}$
if and only if $q_1$ and $q_2$ re in the same orbit.

\remark{Remark} There are also others properties verified, but since
we will not use explicitly here it is unusefull to recall them (cfr.
[Ne] page 138).\endremark

Let $R^{ss}$ be the subset of $q\in R$ such that $\Cal U_q$ is 
a semistable bundle over $X$.

Let $Z_N$ be the variety
$Gr(r;W)\times\dots\times Gr(r;W)$ ($N$ times), where $Gr(r;W)$ is
the grassmannian of quotient vector spaces of $W$ of dimension $r$.
By the tensor product of the Pl\"uker embedding we can find a 
$SL(W)$--equivariant embedding
$Z_N\buildrel {i}\over\hookrightarrow\P(\bigotimes_{i=1}^N\bigwedge^rW)$;
so we can speak about the semistable points under the diagonal action
of $SL(W)$ on $Z_N$.

Let $Z_N^{ss}$ be the open set of $SL(W)$--semistable points 
of $Z_N$ under this action.

If $x\in X$, we can define the map
$$\eqalign{\tau_x :R&\longrightarrow Gr(r;W)\cr
   q&\longrightarrow\Cal U_q\vert_x\cr}$$
where $\Cal U_q\vert_x$ is the fibre at $x$ of the bundle
$\Cal U_q$. The map $\tau_x$is a $SL(W)$--equivariant map.

There is an integer $N=N(r;d;g)$ and $N$ points $x_1;\dots;x_N$ on 
$X$ such that the following properties are verified (cfr. [Ne] page 141):

-- the map
$$\eqalign{\tau_N=\tau_{x_1}\times\dots\times\tau_{x_N}:R&
\longrightarrow Z_N\cr
q&\longrightarrow (\Cal U_q\vert_{x_1};\dots;\Cal U_q\vert_{x_N})\cr}$$
is a $SL(W)$--equivariant map;

-- $R^{ss}=\tau_N^{-1}(Z_N^{ss})$.

\remark{Remark} Again, there are other properties but we will not quote
them here (cfr. [Ne] page 142).\endremark

Now we start the arithmetic construction:

Let $\OK$ be the ring of integers of $K$ and $\Cal W$ be a locally
free model of $W$ over $\OK$. We choose some metrics on $\Cal W$
and a compactification $\T$ of the $GL(\Cal W)$--representation
$\bigotimes_{i=1}^N\bigwedge^r\Cal W$.

Let $Gr(r;\Cal W)$ be the grassmannian of locally free quotients of
rank $r$ of $\Cal W$ and $\Cal Z_N=(Gr(r;\Cal W))^N$. The arithmetic 
scheme $\Cal Z_N$ is a smooth projective model of the variety $Z_N$ over
$\Spec(\OK)$. The group scheme $SL(\Cal W)$ acts on $\Cal Z_N$
and as before we have an $SL(\Cal W)$--equivariant embedding
$$\iota_N:\Cal Z_N\longrightarrow \P(\overline{\Cal W}_T)$$
which extends the $i_N$ defined over $K$.

Let $\overline L=\overline{\Cal O_{\P}(1)}$ be the universal
hermitian line bundle over $\P(\overline{\Cal W_T})$. Let
$\Cal Z_N^{ss}$ be the open set of semistable points of $\Cal Z_N$
and let $\Cal Y$ be the categorical quotient of $\Cal Z_N^{ss}$
(cfr. [Se]). The scheme $\Cal Y$ is a projective scheme over
$\Spec(\OK)$ with a ample line bundle $\Cal L$ such that,
if $\pi:\Cal Z_N^{ss}\to\Cal Y$ is the projection,
$\pi^{\*}(\Cal L)=L^d\vert_{\Cal Z_N^{ss}}$ for some $d>0$.

we have then the following diagram
$$
\CD
R^{ss}@>{\tau_N}>>Z_N^{ss}@>>>\Cal Z_N^{ss}\\
@V {p}VV           @V {\pi_K}VV   @V {\pi}VV\\
\U_X(r;d)@> {\alpha}>>\Cal Y_K@>>>\Cal Y\\
\endCD
$$
Since $\overline L$ is an hermitian line bundle, we can construct a
metric on $\Cal L$ by the formula
$$\Vert m\Vert (x)=\sup_{y\in \pi^{-1}(x)}(\Vert \pi^{\*}(m)\Vert (y))$$
(cfr. [Zh1]).

The line bundle $\Cal M=\alpha^{\*}(\Cal L_K)$ is an ample line
bundle on $\U_X(r;d)$ and we have the formula
$p^{\*}(\Cal M)=\tau_N^{\*}(L^{d})$  (cfr.J [DN]).
So for $[\Cal E]\in \U_X(r;d)$ we can define
$$h([\Cal E])= {{1}\over{d}}h_{\Cal M}([\Cal E])=
{{1}\over{d}}h_{\Cal L}(\alpha ([\Cal E]))=
\inf_{y\in \pi^{-1}(\alpha ([\Cal E]))}
\left\{ h_{\overline L}(y)\right\}.$$
We remark that, $\alpha ([\Cal E])$ is a point in the generic fibre
of the projective scheme $\Cal Y$, so we can speak about its
height as a precise real number (and not just as a number up to
bounded function).

By using Theorem 1 we have then
\proclaim{Proposition 5.3} There exists a universal constant
$C=C(g;r;d)$ such that, for every $[\Cal E]\in\U_X(r;d)$ we
have
$$h([\Cal E])\geq {{rN}\over{[K:\Bbb Q]}}
\cdot {{\widehat{\deg}(\overline{\Cal W})}\over{rk(\Cal W)}}+C.$$
\endproclaim

\

\centerline{\bf References}

\

\item{[A]} S. Ju. Arakelov {\it Intersection Theory of Divisors on a 
Arithmetic Surface},
 Math. U\. S\. S\. R\. Isvestija vol.8, {\bf 6} (1974).

\item{[Bo1]} J.B. Bost, {\it Semistability and heights of cycles}
Inv. Math., {\bf 118} (1994), 223--253.

\item{[Bo2]} J.B. Bost, {\it Heights of stables varieties: Applications
to abelian varieties}, preprint IHES (1994).

\item{[BoGS]} J. Bost, H. Gillet, C. Soul\'e, {\it Heights of
projective varieties and positive Green forms}, Journ. of the AMS,
{\bf 71}, (1994), 903--1027.

\item{[Bu]} J.F. Burnol, {\it Remarques sur la stabilit\'e en
arithm\'etiques}, IMRN, Duke Math. J. {\bf 6}, (1992), 117--127.

\item{[CH]} M. Cornalba, J. Harris, {\it Divisor classes associated
to families of stable varieties with applications to the moduli
space of curves}, Ann. Scient. Ec. Norm. Sup., {\bf 21}, (1988), 455--475.

\item{[DN]} J.-M. Drezet, M.S. Narasimhan, {\it Groupe de Picard des 
Vari\'etes de Modules des Fibr\'es Semi-stables sur les Courbes 
Alg\'ebriques}, Invent. Math.{\bf 97}, (1989). 

\item{[Fa]} G. Faltings, {\it Diophantine approximations on
abelian varieties}, Ann. of Math. {\bf 133} (1991), 549--596.

\item{[Fu]} W. Fulton, {\it Intesection theory}, Springer Verlag, 
(1983)

\item{[FH]} W. Fulton, J. Harris, {\it Representation Theory}, Graduate
Texts in Math. (R.I.M.), {\bf 129}, (1991).

\item{[Ga]} C. Gasbarri, {\it Hauteurs canoniques sur l'espace
de modules des fibr\'es stables sur une courbe algebrique}, Preprint
(1995).

\item{[GS]} H. Gillet, C. Soul\'e, {\it Arithmetic Intersection Theory}, 
Publications Math. IHES, vol. 72, 
(1990).

\item{[Gr]} A. Grothendieck, {\it Techniques de construction en g\'eom\'etrie 
analytique V: Fibr\'es vectoriels, fibr\'es projectifs, fibr\'es
en drapeaux}, S\'em. Cartan 1960/61, No. 12, (1961).

\item{[Ke]} G. R. Kempf, {\it Instability in invariant theory}, Ann.
of Math., {\bf 108}, (1978), 299--316.

\item{[Ma]} V. Maillot, {\it Un calcul de Schubert Arithm\'etique},
Duke Math. J., {\bf 80}, No. 1 (1995), 195--221.

\item{[GIT]} D. Mumford, J. Fogarty, {\it Geometric Invariant Theory}, 2nd 
ed., Springer-Verlag, (1982).

\item{[Mu]} D. Mumford, {\it Stability of projective varieties},
L'Ens. Math., {\bf 23}, (1977), 39--110.

\item{[Ne]} P.E. Newstead, {\it Introduction to moduli problems
and orbit spaces}, Tata Inst. Lect. Notes, Springer Verlag, 1978.

\item{[S]} C. S. Seshadri, {\it Geometric Reductivity over Arbitrary Base},
Advances in Math, {\bf 26}, (1977).

\item{[So]} C. Soul\'e {\it Successive minima on arithmetic varieties},
Comp. Math., {\bf 96}, (1995), 85--98.

\item{[Sz]} L. Szpiro, {\it Degr\'es, Intersections, Hauteurs}, 
S\'eminaire sur les Pinceaux 
Ari\-thm\'e\-ti\-ques:
 la Conjecture de Mordell (L. Szpiro ed.), Expos\'e 1, 
Asterisque 127, Paris, (1985).

\item{[Zh1]} S. Zhang, {\it Geometric reductivity over 
archimedean places}, IMRN, {\bf 10}, (1994), 425--433.

\item{[Zh2]} S. Zhang {\it Heights and reductions of semistable
varieties}, Preprint, Princeton Univ. (1995).

\enddocument